\def\addforarchive{ 
                   \begin{picture}(0,0)
                   \put(260,300){{\small\sf math.QA/0602077}}
                   \put(260,285){{\small\sf KCL-MTH-06-02}}
                   \put(260,270){{\small\sf
                   Hamburger$\;$Beitr\"age$\;$zur$\;$Mathematik$\;$Nr.$\;$231}}
                   \put(260,255){{\small\sf ZMP-HH/05-30}}
                   \end{picture}
                   }

\documentclass{amsproc}

\usepackage{epsfig,amssymb,bbm}

\theoremstyle{definition}

\def\calb          {\ensuremath{\mathcal B}}
\def\calc          {\ensuremath{\mathcal C}}
\def\calh          {\ensuremath{\mathcal H}}
\def\calm          {\ensuremath{\mathcal M}}
\def\cals          {\ensuremath{\mathcal S}}
\def\calv          {\ensuremath{\mathcal V}}
\def\bihom         {bi\-cha\-rac\-ter}
\def\C             {\ensuremath{\mathbbm C}}
\def\CAA           {\ensuremath{\mathcal C_{\!A|A}}}
\def\cft           {conformal field theory}
\def\Chi           {\vartheta}
\def\cir           {\,{\circ}\,}
\def\Cob           {\mbox{\sl Cob}}
\def\Cf            {\mbox{\sl Cor}}
\def\eq            {\,{=}\,}
\def\Hom           {{\rm Hom}}
\def\hy            {$\mbox{-\hspace{-.66 mm}-}$\linebreak[0]}
\def\id            {{\rm id}}
\def\iN            {\,{\in}\,}
\def\io            {\mbox{\sl\i}}
\def\iox           {\ensuremath{\io(\X)}}
\def\M             {\ensuremath{\mathrm M}}
\def\N             {\ensuremath{\mathbb N}}
\def\oti           {\,{\otimes}\,}
\def\Pic           {{\rm Pic}}
\newcommand\quatsch[1]{#1 }
\renewcommand\quatsch[1]{}
\def\R             {\ensuremath{\mathbb R}}
\def\Rep           {{\mathcal R}ep}
\def\Tau           {\ensuremath{\mathcal T}}
\def\top           {\mbox{\sl tft}}
\def\X             {\ensuremath{\mathrm X}}
\def\Y             {\ensuremath{\mathrm E}}
\def\Z             {\ensuremath{\mathbb Z}}
\newcommand\void[1]{}

\theoremstyle{remark}

\begin{document}

\title[Twining characters and Picard groups in RCFT]
{Twining characters and Picard groups \\ in rational conformal field theory}

\author[C.~Schweigert]{Christoph Schweigert}
\address{Organisationseinheit Mathematik, Universit\"at Hamburg,
Bundesstra{\rm\ss{}}e 55, D\,--\,20146 Hamburg}
\email{schweigert@math.uni-hamburg.de}

\author[J.~Fuchs]{J\"urgen Fuchs}
\address{Avdelning fysik, Karlstads Universitet,
Universitetsgatan~5, S\,--\,65188 Karlstad}
\email{jfuchs@fuchs.tekn.kau.se}

\author[I.~Runkel]{Ingo Runkel}
\address{Department of Mathematics, King's College London,
Strand, London WC2R 2LS}
\email{ingo@mth.kcl.ac.uk}

\thanks{C.S.\ is supported by the DFG project SCHW 1162/1-1,
and J.F.\ by VR under project no.\ 621--2003--2385.}

\quatsch{
\dedicatory{This paper is undedicated to the mathematics department
of Hamburg University on the occasion of I.R.'s habilitation.}}

\keywords{Rational conformal field theory, modular tensor categories,
topological field theory, vertex algebras}

\subjclass[2000]{81T40,18D10,18D35,81T45}

\copyrightinfo{2006}{Ribbon Graphs United}

\date{}

\begin{abstract}
Picard groups of tensor categories play an important role in rational conformal
field theory. The Picard group of the representation category $\mathcal C$ of a
rational vertex algebra can be used to construct examples of (symmetric special)
Frobenius algebras in $\mathcal C$. Such an algebra $A$ encodes all data needed 
to ensure the existence of correlators of a local conformal field theory. 
\\
The Picard group of the category of $A$-bimodules has a physical 
interpretation, too: it describes internal symmetries of the conformal field 
theory, and allows one to identify generalized Kramers-Wannier dualities of the
theory.
\\
When applying these general results to concrete models based on affine Lie 
algebras, a detailed knowledge of certain representations of the modular group 
is needed. We discuss a conjecture that relates these representations to those 
furnished by twining characters of affine Lie algebras.
\end{abstract}

\maketitle
\addforarchive

\section{Modular tensor categories and conformal blocks}

Vertex algebras constitute a concise formalization of the
heuristic concept of a chiral symmetry algebra in conformal field theory.
They are, by now, well established as a mathematical structure, and have
applications in various mathematical disciplines. One aspect that makes their 
theory particularly rich is the existence of a notion of tensor product for 
modules over a vertex algebra. The theory of such tensor products was 
developed in a series of articles \cite{hule3.5,hule1,hule2,hule3,huan3}
which form one of Jim Lepowsky's important contributions to the field. 
It closely links vertex algebras to tensor categories.

In recent work \cite{huan24}, Huang extends this relationship by identifying a 
class of vertex algebras for which the representation category carries the 
stronger structure of a modular tensor category.
This justifies the point of view, initiated in the work of Moore and Seiberg
\cite{mose3}, that modular tensor categories furnish an axiomatization of
the chiral data of rational two-dimensional conformal quantum field theory.

For the present contribution, we adopt the following definition of a modular 
tensor category: it is an abelian, semi-simple, $\C$-linear monoidal category 
$\calc$, equipped with a functorial braiding 
$\,c_{X,Y}{:}\ X\oti Y \,{\xrightarrow{\sim}}\, Y \oti X$ and functorial
twist isomorphisms $\theta_X{:}\ X \,{\xrightarrow{\sim}}\, X $, subject
to the following axioms:
   \def\leftmargini{1.1em}~\\[-1.2em]
   \begin{itemize}\addtolength\itemsep{2pt}
\item Braiding and twist are compatible in the sense that $\theta_{X\otimes Y}
      \eq c_{Y,X}^{}{\circ}\,(\theta_Y{\otimes}\theta_X)\cir c_{X,Y}^{}\!$.
\item There exists a compatible duality $X\,{\mapsto}\, X^\vee$ on $\calc$.
\item The tensor unit {\bf1} is simple.
\item There are only finitely many isomorphism classes of simple objects.
 \\[1pt] We denote a set of representatives for these classes by $\{U_i\}_{i\in I}$.
\item Finally, there is a non-degeneracy condition on the braiding.
 \\[1pt] It ensures in particular that the natural transformations of the
      identity functor on $\calc$ are controlled by the fusion ring
      $K_0(\calc)$, in the sense that 
$${\rm End}(\mbox{\sl Id}_{\calc}) \cong K_0(\calc)\otimes_{\Z}^{} \C \,. $$
\end{itemize}

By results of Reshetikhin and Turaev \cite{retu}, to each modular tensor 
category one can associate a monoidal category $\Cob_\calc$ of decorated
cobordisms and a monoidal functor $\top_\calc$ from $\Cob_\calc$ (or rather, a 
version of this category that is enriched over $\mathcal Vect_\C$) into the 
category $\mathcal Vect_\C$ of fi\-ni\-te-di\-men\-si\-o\-nal complex vector 
spaces. The objects of $\Cob_\calc$ are called extended surfaces; an extended 
surface \Y\ is an oriented closed two-dimensional topological manifold, together
with a Lagrangian subspace of $H_1(\Y,\R)$ and germs of disjoint arcs that are 
labeled by objects of $\,\calc$. The finite-dimensional vector spaces $\calh(\Y)
\eq \top_\calc(\Y)$ are called spaces of {\em conformal blocks\/}. As we will 
see, these vector spaces are expected to capture important properties of the 
(sheaves of) conformal blocks constructed from vertex algebras.  

The morphisms of $\Cob_\calc$ are cobordisms $(\M,\partial_-\M,\partial_+\M)$ 
containing a ribbon graph (``Wilson lines'', in physicists' language) which is 
required to end at the germs of arcs on the boundary $\partial\M\eq
\partial_-\M\,{\sqcup}\,\partial_+\M$. Segments of the ribbon graph are labeled 
by objects of $\calc$, and vertices by suitable morphisms of $\calc$. 
It follows immediately from the properties of a monoidal functor that for a 
cobordism of the form $(\M,\emptyset,\partial\M)$ the functor $\top_\calc$
gives a linear map
   $$ \top_\calc(\M) \equiv \top_\calc(\M,\emptyset,\partial \M)
   : \quad \C \to \calh(\partial\M) \,. $$
Thus this type of 
cobordism allows us to specify vectors in the space of conformal blocks 
associated to the boundary of a three-manifold with embedded ribbon graph.

Apart from monoidality, the functor $\top_\calc$ must obey a few more axioms 
(see e.g.\ \cite{fffs3} for the explicit statement) -- naturality, 
multiplicativity, normalization, and functoriality, i.e.\ compatibility with 
gluing: if the cobordism $\M$ is obtained by gluing
two cobordisms $\M_1,\M_2$ with a morphism 
$$ 
f:\ \ \partial_+\M_1\to \partial_-\M_2 
$$
of extended surfaces, then the corresponding linear maps $\top_\calc(\cdot)$ 
are related by
$$ 
\top_\calc(\M) = \kappa^m \, 
\top_\calc(\M_2) \circ f_\sharp \circ \top_\calc(\M_1) 
$$
with some complex number $\kappa$ that depends only on $\calc$ and an integer 
$m$ that depends on the cobordisms $\M_1$ and $\M_2$, and on $f$.
The functor $\top_\calc$ 
thus provides a three-dimensional topological field theory \cite{TUra}.

The axioms for $\top_\calc$ have two important consequences: 
   \def\leftmargini{1.1em}~\\[-1.2em]
   \begin{itemize}\addtolength\itemsep{5pt}
\item They endow $\calh(\Y)$ with a projective representation of the
      mapping class group ${\rm Map}(\Y)$ of \Y.
\item Taking a cylinder $\Y \,{\times}\, [0,1]$ over an extended surface \Y\
      and gluing two of its arcs
      by a tube on one side of the resulting cobordism allows one to obtain
      factorization isomorphisms for conformal blocks.
\end{itemize}

The term ``conformal block'' is also used in a mathematically a priori 
different context. Given a conformal vertex algebra $\calv$ and an $m$-tuple 
$(\lambda_1,\lambda_2,...\, ,\lambda_m)$ of $\calv$-modules, one 
can construct a vector bundle $\calb(\lambda_1,\lambda_2,...\, ,\lambda_m)$ 
over the moduli space $\calm_{g,m}$ of complex curves of genus $g$ with $m$ 
marked points. The action of the Virasoro algebra on these representations 
leads to a projectively flat connection on $\calb$, the Knizhnik-Zamolodchikov 
connection. On the other hand, the structural morphisms that, for a sufficiently
nice subclass of vertex algebras \cite{huan24}, endow the representation 
category $\Rep(\calv)$ with the structure of a modular tensor category are 
deduced solely from monodromy data of conformal blocks on $\C{\mathbb P}^1$. 

Thus a rational vertex algebra, i.e.\ a vertex algebra whose representation 
category is a modular tensor category $\calc$, supplies us with two a priori
unrelated representations of the mapping class group $\pi_1(\calm_{g,m})$: 
the one obtained from the mo\-no\-dromies of the vector bundles $\calb$ over 
$\calm_{g,m}$, and the one provided by the functor $\top_\calc$. We will, from 
now on, assume that these two representations coincide for all genera $g$ and 
all values of $m$. This assumption implies in particular an explicit formula for
the rank of the vector bundles $\calb(\lambda_1,\lambda_2,\ldots,\lambda_m)$, 
the Verlinde formula, which is of interest to algebraic geometers (see 
\cite{sorg} for a review). Factorization rules of conformal blocks provide 
constraints on the representations of mapping class groups. Our assumption is 
thus closely linked to factorization rules; the latter have been established 
for specific classes of models, see e.g.\ \cite{tsuy,FRbe}.

The horizontal sections of $\calb$ are, in general, multivalued. Therefore 
their relation to correlation {\em functions\/} of a conformal field theory 
calls for a clarification. Our assumptions about the monodromies of $\calb$ 
allow us to address this question using the functor $\top_\calc$.

\section{Some geometry for local conformal field theories}

Concretely, the desired clarification is provided by the following geometric 
construction. We are interested in a local CFT on a two-dimensional manifold
\X\ that can have a boundary. One must actually carefully distinguish between
two different types of local CFTs: one in which the surface \X\ is required
to be oriented, and one in which no orientation of \X\ is presupposed and
accordingly also unorientable surfaces are admitted. However, for our
purposes these two types of local CFTs can be treated simultaneously.

The first step is to replace \X\ by the pair $(\widehat\X,\sigma)$ consisting
of the complex double $\widehat\X$ of \X\ and an orientation reversing
involution $\sigma$. Notice that \X\ can be recovered as the quotient
by the action of the group $\langle \sigma\rangle$: there is a projection
$$  
\pi: \ \ \widehat\X \,\mapsto\,\X\,\cong\,\widehat\X\,/\langle\sigma\rangle \,. 
$$

We can now state the `principle of holomorphic factorization' which relates 
local CFT on \X\ to chiral CFT on $\widehat\X$ (for certain
classes of conformal field theories this principle can formally be derived from
an action functional \cite{witt39}): a correlator $\Cf(\X)$ on \X\ is a specific 
vector in $\calh(\widehat\X)$. These vectors must obey two additional axioms:
   \def\leftmargini{1.1em}~\\[-1.2em]
   \begin{itemize}\addtolength\itemsep{4pt}
\item They must be invariant under the action of $\,{\rm Map}(\X)\,{\cong}\,
      {\rm Map}(\widehat \X)^\sigma_{\phantom|}$.
   \\ This group, also called the relative modular group \cite{bisa2}, acts
      genuinely on $\calh(\widehat\X)$, rather than only projectively.
\item They must possess certain factorization properties.
   \\[1pt] We refer to
      \cite{fjfrs,fjfrs2} for a precise formulation of these constraints.
\end{itemize}

In order to apply the relation to TFT outlined in section 1, we
restrict our attention to {\em rational\/} full conformal field theories
(RCFTs), that is, full CFTs whose correlators are obtained from
the conformal blocks of a rational vertex algebra.
     
To put holomorphic factorization to work we now look for a cobordism
$(\M_\X,\emptyset,\widehat \X)$ such that the vector
$\top_\calc(\M_\X,\emptyset,\widehat\X)1 \iN \calh(\widehat\X)$
is the correlator $\Cf(\X)$. It turns out that the following quotient of
the interval bundle on $\widehat\X$ is appropriate:
$$  
\M_\X = \big(\,\widehat\X \,{\times}\, [-1,1]\,\big) \,{\slash}\,
\langle(\sigma,t\,{\mapsto}\, {-}t)\rangle  \,. 
$$
This three-manifold is oriented, has boundary $\partial\M_\X\,{\cong}\,
\widehat\X$, and it contains \X\ as a retract: the embedding $\io$ of \X\ is 
to the fiber $t{=}0$, the retracting map contracts along the intervals.

The details of the construction are summarized in appendix A of \cite{fjfrs}. 
Here we content ourselves with emphasizing that it involves placing a 
uni-trivalent ribbon graph $\Gamma_{\!\X}$ on the union of $\iox$ with the 
intervals over boundary points of $\iox$. This graph contains in particular 
$\io(\partial\X)$. Univalent vertices of $\Gamma_{\!\X}$ are required to lie 
on $\pi^{-1}(\partial \X) \,{\subset}\, \widehat\X \,{\cong}\,\partial\M_\X$.

\section{Frobenius algebras}

The labeling of the ribbon graph $\Gamma_{\!\X}$ requires further data. The 
central idea of the TFT approach to RCFT correlators, as developed in
\cite{fuRs4,fuRs8,fuRs10,fjfrs} is that these data are provided
by a (symmetric special) Frobenius algebra $A$ in the tensor category $\calc$. 

A {\em Frobenius\/} algebra $A$ in $\calc$ is, by definition, an object of 
$\calc$ carrying the structures of a unital associative algebra and of a 
counital coassociative coalgebra in $\calc$, satisfying the compatibility 
requirement that the coproduct $\Delta{:}\ A \,{\to}\, A\oti A$ is a morphism of 
$A$-bimodules (or, equivalently, that the product $m{:}\ A\oti A\,{\to}\,A$ is 
a morphism of $A$-bi-comodules). A Frobenius algebra is called {\em special\/} 
iff the coproduct is a right-in\-ver\-se to the product -- this means in 
particular that the algebra is separable -- and a multiple of the unit is a 
right-inverse to the counit. $A$ is called {\em symmetric\/} iff the two 
isomorphisms $A\,{\to}\, A^\vee$
that are naturally induced by product, counit and duality coincide.

We mention that in case \X\ is unoriented, still more structure is needed:
$A$ must then be a {\em Jandl\/} algebra, i.e.\ a symmetric special Frobenius
algebra coming with an algebra isomorphism $A \,{\to}\, A^{{\rm opp}}$ that 
squares to the twist $\theta_{\!A}$. This turns out to be the appropriate 
generalization 
of the notion of an algebra with involution to braided tensor categories.

Given such a symmetric special Frobenius algebra $A$, one can identify
boundary conditions of the CFT with $A$-modules. Boundary fields $\Psi^{MN}_U$
that change the boundary condition from $M$ to $N$ and have chiral label
$U$ are in bijection with $A$-module morphisms in $\Hom_A(M\oti U, N)$. 
Types of (topological) defect lines can be identified with isomorphism 
classes of $A$-bimodules. The disorder fields $\Psi^{B_1B_2}_{UV}$ which change
a defect line from type $B_1$ to type $B_2$ carry two chiral labels $U$ and $V$.
They are in bijection with $A$-bimodule morphisms in
$\Hom_{A|A}( U\,{\otimes^+}\,B_1\,{\otimes^-}\,V, B_2 )$. Here the 
left and right $A$-actions $\rho_{l/r}$ on the bimodule $B_1$ are used to define 
a bimodule structure on the object $U\oti B_1 \oti V$ by 
$(\id_U\oti\rho_l\oti\id_V)\cir (c^{-1}_{U,A}\oti\id_{B_1}\oti\id_V)$ and 
$(\id_U\oti\rho_r\oti\id_V)\cir (\id_U\oti\id_{B_1}\oti c^{-1}_{\!A,V})$, 
respectively. 

Notice that $A$ itself has a canonical structure of an $A$-bimodule; in fact, 
it is the tensor unit in the tensor category \CAA\ of $A$-bimodules.
In applications to conformal field theory, $A$ is required to be simple
as a bimodule over itself. In the CFT context, we also call $A$ the
invisible defect. Bulk fields of the CFT are the special disorder fields 
that correspond to $B_1\eq B_2\eq A$. As a consequence, the degeneracy of 
bulk fields with chiral labels $U_i$ and $U_j$ is given by
\begin{equation}
Z_{ij} := \dim_\C^{} \Hom_{A|A}(U_i\,{\otimes^+}A\,{\otimes^-}\,U_j, A) \,.
\label{pf} 
\end{equation}
Thus such expressions yield modular invariant bulk field partition functions
of CFTs, a result that has been obtained earlier \cite{boek3} in the 
framework of subfactor theory.

The Frobenius algebra $A$ is, of course, also a module over itself and thus
corresponds to a boundary condition. But in fact any other $A$-module $M$
supplies us with a symmetric special Frobenius algebra as well, namely the
internal End $\underline{{\rm End}}(M)$. All these algebras 
$\underline{{\rm End}}(M)$ are Morita equivalent to $A$, and in the TFT 
approach they all provide equivalent descriptions of a full local CFT.

\smallskip

On the basis of these informations one can establish, for instance, the 
following further results:
   \def\leftmargini{1.6em}~\\[-1.2em]
   \begin{itemize}\addtolength\itemsep{3pt}
\item[I\ \ ] One can derive concrete expressions for partition functions
     of boundary, bulk and defect fields. Their coefficients are non-negative
     integers, and they satisfy consistency requirements like being modular
     invariant and forming NIMreps of the fusion rules \cite{fuRs4}. 
  \\ (However, not every bilinear combination of characters that is invariant
     under the action of the modular group is indeed the partition function
     of a conformal field theory. Many counterexamples are known.)
\item[II\ ] These results can be extended to unorientable surfaces,
       for which the Frobenius algebra must be a Jandl algebra 
       \cite{fuRs8,fuRs11}.
\item[III] The expressions for correlation functions can be made particularly
      explicit \cite{fuRs9} for 
     theories with torus partition function of simple current type;
     some details will be given below.
  \\ It is worth stressing, though, that the TFT approach to RCFT correlators
     treats all RCFTs -- the simple current case and those with 
     exceptional modular invariants -- on an equal footing.
\item[IV] One can derive explicit expressions for the
     coefficients of operator product expansions \cite{fuRs10}.
\item[V\ ] Finally, one can prove \cite{fjfrs} modular invariance, at arbitrary 
     genus, and factorization of the correlators obtained in the TFT approach.
\end{itemize}

\section{Picard groups}

{}From the discussion above we learn in particular that in any full local 
conformal field theory two tensor categories are encountered naturally: the 
(modular) tensor category $\calc$ that describes the chiral data, and the tensor 
category \CAA\ of $A$-bimodules that describes defects and their fusion. 
By standard arguments, the latter does not depend on the choice of $A$
within one and the same Morita class. It should be noted that \CAA\ 
is not, in general, a braided tensor category. This fits well with the absence
of any physically reasonable notion of a braiding of defects.
The following construction can be applied to both tensor categories $\calc$ 
and \CAA.

An object $V$ of a tensor category with simple tensor unit is called 
{\em invertible\/} iff there exists an object $W$ such that $V\oti W\,{\cong}\, 
{\mathbf 1}\,{\cong}\, W \oti V$. The tensor product of $\calc$ endows
the set of isomorphism classes of invertible objects with the structure of
a group, the Picard group $\Pic(\calc)\,{\subset}\,K_0(\calc)$. The Picard group
of a braided tensor category is abelian. If $\calc$ is the representation 
category of a rational vertex algebra, the term `simple currents' 
\cite{scya6} is also used for elements of the Picard group.

The following two results allow one to control aspects of tensor categories,
given their Picard groups:
   \def\leftmargini{1.1em}~\\[-1.2em]
   \begin{itemize}\addtolength\itemsep{3pt}
\item Let $G$ be a finite group. Monoidal categories $\calc$ that are enriched
      over $\mathcal Vect_\C$ 
      together with an isomorphism $K_0(\calc)\,{\to}\, \Z G$ are in
      bijection with elements of the cohomology group $H^3(G,\C^\times)$.
\item Let $G$ be a finite abelian group. Braided monoidal categories $\calc$ 
      enriched over $\mathcal Vect_\C$ 
      together with an isomorphism $K_0(\calc)\,{\to}\,\Z G$ are in
      in bijection with elements of the group $H^3_{{\rm ab}}(G,\C^\times)$
      of Eilenberg-Mac Lane's \cite{eima2} abelian group cohomology.
\end{itemize}

It turns out that elements of $H^3_{{\rm ab}}(G,\C^\times)$ are in bijection
with quadratic forms on $G$. In the case at hand, the relevant quadratic form
on the Picard group is given by the eigenvalue of the twist, i.e.\ by (the 
exponential of minus $2\pi{\rm i}$ times) the fractional part of the 
conformal weight.

It is now time to present some examples for symmetric special Frobenius 
algebras. The tensor unit is such an algebra and, more generally, for any object
$U$ of $\calc$ with $\dim(U)\,{\ne}\,0$ the object $U \oti U^\vee$ can be 
endowed with such an algebra structure. These algebras are all Morita 
equivalent; they yield in (\ref{pf}) the (charge-) diagonal bulk field partition
function, $Z_{ij} \eq\delta_{i,j^\vee_{}}$. This situation has been termed the 
``Cardy case'' in the physics literature.

The next tractable class of symmetric special Frobenius algebras are 
{\em Schel\-le\-kens\/} algebras. A Schellekens algebra is a symmetric special 
Frobenius algebra for which each simple subobject is invertible and which is 
simple as a left-module over itself. Every such algebra can be constructed 
from a subgroup $H\,{\leq}\, \Pic(\calc)$ (to be called the support of $A$) 
and a cochain $\omega{:}\ H\,{\times}\, H\,{\to}\,\C^\times$ with
${\rm d}\omega \eq \psi\raisebox{-2pt}{$|{}_H$}$, where $\psi$ represents the 
class in $H^3(G,\C^\times)$ that determines the monoidal structure of the 
Picard subcategory of $\calc$, i.e.\ of the full tensor subcategory of $\calc$ 
generated by the invertible objects. For any symmetric special Frobenius 
algebra that is simple as a left module over itself, one has the estimate
\cite{gann17,fuRs9}
$$ 
\dim_{\C} \Hom(U,A) \leq \dim(U) \,; 
$$
as a consequence, the object underlying a Schellekens algebra is a direct sum 
$\bigoplus_{h\in H} U_h$ of pairwise non-isomorphic invertible objects.

It is, in fact, appropriate to regard Schellekens algebras as the generalization
to braided categories of twisted group algebras. The latter are classified
by $H^2(G,\C^\times)$; for an abelian group $G$, this cohomology
group is isomorphic to the group $A\!B(G,\C^\times)$ of alternating \bihom s of
$G$, an isomorphism being given by
$$
\begin{array}{rcl}
H^2(G,\C^\times) &\!\!\to\!\!& A\!B(G,\C^\times) \\{}\\[-10pt]
{} [\omega] &\!\!\mapsto\!\!& \xi(g,h)
:=\mbox{\Large$\frac{\omega(g,h)}{\omega(h,g)}$} \,.  \end{array} 
$$
While an alternating \bihom\ $\xi$ is characterized by $\xi(g,g)\eq1$, its 
appropriate braided generalization, a {\em Kreuzer-Schellekens bihomomorphism\/}
(KSB), is a \bihom\ $\Xi$ whose diagonal values are given by the eigenvalues 
of the twist on the invertible objects, $\Xi(g,g) \eq \theta_g$.

Next we realize that the multiplication on a Schellekens algebra provides us 
with a KSB. Indeed, for a Schellekens algebra of support $H\,{\le}\,
\Pic(\calc)$, pick for each group element $[U_g]\iN H$ a nonzero morphism 
$\iota_g\iN\Hom(U_g,A)$. The space $\Hom(U_g\oti U_h, U_{gh})$ is 
one-dimensional, and hence there exist numbers $\Xi_A (h,g)$ such that 
$$ 
m \cir c_{A,A}^{} \cir (\iota_g\oti\iota_h) = \Xi_A (h,g)\,
m \cir (\iota_g\oti \iota_h) \,. 
$$
It should be appreciated that the resulting function
$$ 
\Xi_A :\ \ H\,{\times}\, H \to \C 
$$ 
does not depend on the choice of the morphisms $\iota_g$, but depends on the 
multiplication on $A$. One shows that $\Xi_A (h,g)$ is a KSB, and that a 
Schellekens algebra is uniquely characterized, up to isomorphism, by the 
support $H$ of $A$ and the choice of a KSB on $H$.

We are now in a position to express quantities of interest in conformal field
theory in terms of the KSB. First, we obtain the partition function 
\begin{equation}
Z_{ij}(A) = \frac1{|H(A)|}\! \sum_{\ g,h\,\in H(A)}\!\!\!
\Chi_{i}^{}(h) \,\, \Xi_A(h,g) \, \delta_{j^\vee_{}\!,gi}  
\label{Z} \end{equation}
for bulk fields, which reproduces a result of \cite{krSc}.
Here $\Chi_{i}^{}(h)$ is the so-called (exponentiated) monodromy charge of
$U_i$ with respect to $h$, which is the combination
$$ 
\Chi_{i}^{}(h) := \theta_{hi}^{}\,\theta_h^{-1}\,\theta_i^{-1} 
$$
of twist eigenvalues.

Boundary conditions can be determined from the decomposition of induced 
$A$-mo\-dules into simple $A$-modules. They correspond to orbits of 
the support $H(A)$ on (isomorphism classes of) simple objects of $\calc$ and 
representations of the twisted group algebra $\C_{\epsilon_U^{}}\cals_U$, where 
$\cals_U\,{\leq}\,H(A)$ is the stabilizer of the simple object $U$ and the 
alternating \bihom\ $\epsilon_U$ on $\cals_U$ is given by
$$
\epsilon_U(g,h) = \phi_U(g,h) \, \Xi_A(h,g) \,,
$$
i.e.\ is the product of the restriction of the KSB $\Xi_A$ to $\cals_U$ and a 
certain two-cochain $\phi_U$. The latter is a gauge independent $6j$-symbol, 
defined by
\begin{equation}
\gamma_g \cir (\id_{U_g}^{} \oti \delta_h)
=: \phi_U(g,h) \, \delta_h \cir (\gamma_g \oti \id_{U_h}) \,,
\label{P}  \end{equation}
with $\{\gamma_g\}$ and $\{\delta_h\}$ (arbitrary) bases of the
one-dimensional morphism spaces $\Hom(U_g\oti U,U)$ and $\Hom(U\oti U_h,U)$,
respectively. A similar description exists for defect lines. Also, explicit 
formulae for boundary states have been derived in \cite{fuRs9},
which prove conjectures made in \cite{fhssw}.

\section{Twining characters}

The gauge independent $6j$-symbols $\phi_U$ defined in formula \eqref P deserve 
further study. One can show that, for a general modular tensor category, they 
are \bihom s, in fact KSBs, on the stabilizer $\cals_U$. To proceed, consider 
the special case of chiral WZW conformal field theories: As is well known, to 
each pair $({\mathfrak g},k)$ with $\mathfrak g$ a finite-dimensional complex
simple Lie algebra and $k\iN \N$, one can associate a modular tensor category 
$\calc({\mathfrak g},k)$. It can be realized as the category of integrable
highest weight modules at level $k$ of the untwisted affine Lie algebra
$\hat {\mathfrak g}$ with horizontal subalgebra ${\mathfrak g}$, or as
the representation category of the vertex algebra $\calv({\mathfrak g},k)$ 
that can be defined on the irreducible highest weight 
${\mathfrak g}$-module with highest weight $k\Lambda_{(0)}$.

For these categories, conformal field theory leads to a conjecture, to be
presented as relation \eqref{con.phi} below, that links the $6j$-symbols 
$\phi_U(g,h)$ to representation theoretic quantities, the so-called twining 
characters \cite{fusS3,furs} (or graded traces, see \cite{dolm13}).
A symmetry $\dot\omega$ of the Dynkin diagram of a Kac-Moody algebra (or 
generalized Kac-Moody algebra) $\hat{\mathfrak g}$ gives rise to an outer 
automorphism $\omega$ of $\hat{\mathfrak g}$ that respects a triangular 
decomposition. We call such an automorphism of $\hat{\mathfrak g}$ a diagram 
automorphism.  It should be noted that in the case that $\hat{\mathfrak g}$ 
is an untwisted affine Lie algebra and $\dot\omega$ acts nontrivially on the 
node corresponding to the vacuum representation, these Lie algebra 
automorphisms do {\em not\/} give rise 
to automorphisms of the vertex algebra $\calv({\mathfrak g},k)$.

Now recall that the simple objects of $\calc({\mathfrak g},k)$ are irreducible 
highest weight representations of the untwisted affine Lie algebra
$\hat{\mathfrak g}$. As for any tensor category, the Picard group of 
$\calc({\mathfrak g},k)$ acts on the isomorphism classes of simple objects by 
the tensor product. Since an invertible object $L$ obeys $L \oti L^\vee 
\,{\cong}\, {\bf1}$ and since $\dim(L)\eq\dim(L^\vee)$, the (quantum) dimension
 of $L$ fulfils $(\dim(L))^2 \eq 1$. When scanning integrable 
$\hat{\mathfrak g}$-modules with this property one finds \cite{jf15} 
that for any pair $(\mathfrak g,k)$ (with the exception of $(E_8,2)$, which
is to be excluded from the discussion below), every element of the Picard 
group corresponds to a diagram automorphism of $\hat{\mathfrak g}$. 

Given an automorphism $\varpi$ of $\hat{\mathfrak g}$ that preserves a fixed 
Borel subalgebra and an irreducible highest weight representation $(V, \rho)$
with $\rho{:}\ \hat{\mathfrak g}\,{\to}\,\mathfrak{gl}(V)$, another irreducible
highest weight representation $(V,\rho^\varpi)$ is obtained by setting $\rho^
\varpi\,{:=}\, \rho\cir\varpi$. If $\varpi\eq\omega$ is a diagram automorphism, 
then the highest weight spaces of these two representations coincide. 
If the irreducible representations $(V, \rho)$ and $(V,\rho^\omega)$ are 
isomorphic, the representation is called a fixed point of the automorphism
$\omega$. In this case, there exists a unique linear map $\Tau_\omega{:}\ 
V\,{\to}\,V$ that restricts to the identity on the highest weight space 
and obeys $\rho\cir\omega(x)\, \Tau_\omega \eq \Tau_\omega\, \rho(x)$ for all
$x\iN\hat{\mathfrak g}$. Given an irreducible highest weight representation 
$(V_\Lambda, \rho_\Lambda)$ of highest weight $\Lambda$ that is a fixed point 
for some diagram automorphism $\omega$, one defines the 
(Viraso\-ro-spe\-ci\-alized)
{\em twining character\/} of $(V_\Lambda, \rho_\Lambda)$ as the function
\begin{equation}
\raisebox{.15em}{$\chi$}^\omega_\Lambda(\tau) := {\rm Tr}_{V_\Lambda}\,
\Tau_\omega\, \exp\big(2\pi{\rm i}\tau (L_0{-}\mbox{$\frac c{24}$})\big)  
\end{equation}
on the complex upper half plane, where $L_0$ is the action of the zero mode 
of the Virasoro algebra obtained by the Sugawara construction. 

In particular, given an element $g$ of the Picard group $\Pic
(\calc({\mathfrak g},k))$, we obtain a twining character for each integrable 
highest $\Lambda$ that is a fixed point under the corresponding Lie algebra 
automorphism $\omega_g$.  Explicit character formulae for the twining characters 
of irreducible highest weight modules have been derived in \cite{fusS3,furs}. 
(For twining characters of other types, see \cite{kakw,wend2,nait11,nait13}.)
{}From these character formulae, it follows that, under the usual action of 
the modular group $SL(2,\Z)$ on the upper half plane, these twining 
characters transform into linear combinations of twining characters for the 
same automorphism $\omega_g$ -- a non-trivial empirical observation from 
the point of view of Lie theory. In particular one finds
\begin{equation}
\raisebox{.15em}{$\chi$}^{\omega_g}_\Lambda(-\mbox{\large$\frac1\tau$}) 
= \sum_{\Lambda'} S^{\omega_g}_{\Lambda,\Lambda'}\, 
\raisebox{.15em}{$\chi$}^{\omega_g}_{\Lambda'}(\tau) 
\end{equation}
with a unitary and symmetric matrix $S^{\omega_g}_{\phantom:}$.

The following conjecture relates the matrix $S^{\omega_g}_{\phantom:}$ obtained
this way to the gauge independent $6j$-symbols $\phi_U$ defined by \eqref P: 
\begin{equation}
S^{\omega_g}_{\Lambda,h\Lambda'} = \Chi_{\!\Lambda}^{}(h)\, 
\phi_{\Lambda'}^{}(g,h)^{-1} S^{\omega_g}_{\Lambda,\Lambda'} \,, \label{con.phi}
\end{equation}
where $\Chi_{\!\Lambda}^{}(g)^{}\eq \theta_{g\Lambda}^{}\theta_g^{-1}\theta
_\Lambda^{-1}$
is the exponentiated monodromy charge, which we already encountered in \eqref Z. 

A proof of this conjecture seems to be out of reach at present. However, the 
formula has been checked in a huge number of non-trivial cases; most of these 
checks have been performed with a sophisticated computer package \cite{kac}.
For the role of the $6j$-symbols played in the untwisted stabilizer, we also
refer to \cite{bant7,bant6}. The matrices $S^{\omega_g}$ are also an important
ingredient in refinements of the Verlinde conjecture for non-simply connected
Lie groups \cite{fusS6,fuSc8}.

\section{The bimodule Picard group} 

Let us now turn our attention to the Picard group of the category \CAA\ of 
$A$-bi\-modules. This group is of considerable interest, too: it turns out 
\cite{ffrs3} that it describes internal symmetries of a local CFT. Since the 
tensor category \CAA\ is, in general, not braided, this group of symmetries 
can be non-abelian. To give some examples, for the Ising model one has 
$\Pic(\CAA)\,{\cong}\,\Pic(\calc)\,{\cong}\,\Z_2$, while for the critical 
three-state Potts model we find a non-abelian symmetry group, $\Pic(\CAA)\,
{\cong}\,S_3$. These results conincide with the symmetries of the Ising and 
the three-state Potts model on a lattice which, as their names suggest, 
produce the corresponding conformal field theories in the continuum limit.

The action of an invertible bimodule $B$ on a boundary condition $M$ of
the CFT is given by
$$ 
M \mapsto B\,{\otimes_{\!A}}\, M \,, 
$$
where $\otimes_{\!A}$ is the tensor product over $A$.
For the action on boundary fields, we find 
$$ 
\Hom(M_1\oti U, M_2) \,\mapsto\, \Hom((B\,{\otimes_{\!A}}\,M_1)\oti U, 
B\,{\otimes_{\!A}}\,M_2) \,. 
$$
For the action on bulk fields, which is somewhat more involved, we refer to 
\cite{ffrs3,ffrs4}. 

\medskip

If $A \eq A(H,\Xi_A)$ is a Schellekens algebra, then there is a subgroup
of internal symmetries of the form 
$$
H^*_{} \,{\times_H^{}}\,\Pic(\calc) \,, 
$$ 
where $H$ is mapped to its character group $H^*_{}$ by the homomorphism
$H \,{\to}\, H^*_{}\!$, $H\,{\ni}\,h \,{\mapsto}\, \Xi_A(\,\cdot,h)$ and then 
acts by right-multiplication.
This statement can be derived as a special case of the following result:
Suppose that $\calc'$ is a braided tensor (in fact, semisimple ribbon) category 
with the property that $\Pic(\calc')$ acts freely on the isomorphism classes of 
simple objects of $\calc'$. (For any tensor category, the full subcategory of 
all invertible objects has this property.) Let $A$ be a Schellekens algebra in 
$\calc'$ with support $H$. Then the fusion ring of $A$-bimodules is given by
\begin{equation} 
K_0(\calc'_{\!A|A}) \cong K_0(\calc') \otimes_{\Z H}^{} \Z H^*_{} . 
\label{res}\end{equation}
Here the group ring $\Z H$ is regarded as a subring of $K_0(\calc')$, 
and its action on $\Z H^*_{}$ is again defined using the KSB $\Xi_A$ of $A$.

The result (\ref{res}) can be seen as follows (for details see \cite{ffrs5}). 
One first shows that the automorphism group of a Schellekens algebra
$A\eq\bigoplus_{h\in H} U_h$ is canonically isomorphic to the character group
$H^*_{}$. This extends results about the relationship between simple current 
extensions and orbifolds by abelian groups. Next one 
shows that every simple $A$-bi\-mo\-dule is isomorphic to a bimodule of the 
form $A\oti U$ with left action $m\oti\id_U$ and right action $(m\oti\id_U)\cir
(\id_A\oti c_{UA})\cir(\id_{A\otimes U}\oti\psi)$ for some algebra automorphism
$\psi$ of $A$. The tensor product over $\Z H$ in (\ref{res}) takes care 
of an over-parametrization in a way that is compatible with the fusion rules.
%

\medskip
 
These considerations can be generalized to include not only internal 
symmetries, but also duality symmetries of Kramers-Wannier type. The fact that 
$A$ is a Frobenius algebra implies that the monoidal category \CAA\ 
comes with a duality; it is even a sovereign tensor category. The conformal 
field theory has Kramers-Wan\-nier like dualities whenever there is an object
$B\iN\CAA$ such that $B \,{\otimes_{\!A}}\, B^\vee$ is a direct sum of
invertible bimodules, i.e.\ is in the Picard subcategory of \CAA. The 
well-known Kramers-Wannier duality for the Ising model can thus be deduced from 
the fusion rule $\sigma\oti\sigma \,{\cong}\,{\bf 1} \,{\oplus}\, \epsilon$ 
for the primary field $\sigma$ of conformal weight $\frac1{16}$, in agreement 
with another approach to such dualities based on the symmetries of boundary 
states \cite{ruel'5}.

\section{Conclusions}

The TFT approach to the construction of CFT correlators provides a powerful 
algebraization of many questions arising in conformal field theory. 
This allows both for making rigorous statements about rational conformal field 
theories and for setting up efficient algorithms. A rich dictionary that relates
algebraic concepts and physical notions is emerging; the relation between Picard
groups and symmetry groups is merely one example. 

Interesting algebraic structures also show up when a route analogous to the 
TFT approach is taken in the framework of vertex algebras. The results of 
\cite{huko,huko2} should be seen as just the beginning of a rich extension of 
the theory of vertex algebras.

To summarize, a central idea of the construction presented in this article is 
to represent conformal field theory quantities as invariants of knots and 
links in three-manifolds. Tools for achieving this are provided by topological
field theory, which in turn is based on the notion of a modular tensor category.
It seems thus fair to say that the tensor structure present 
for the modules over the underlying vertex algebra, which is due in particular 
to the work of Jim Lepowsky, is a key ingredient in the construction.

\bigskip

 \newcommand\wb{\,\linebreak[0]} \def\wB {$\,$\wb}
 \newcommand\Bi[2]    {\bibitem[#2]{#1}}
 \newcommand\JO[6]    {{\em #6}, {#1} {#2} ({#3}), {#4--#5} }
 \newcommand\J[7]     {{\em #7}, {#1} {#2} ({#3}), {#4--#5} {{\tt [#6]}}}
 \newcommand\JJ[6]    {{\em #6}, {#1} {#2} ({#3}), {#4} {{\tt [#5]}}}
 \newcommand\Pret[2]  {{\em #2}, pre\-print {\tt #1}}
 \newcommand\BOOK[4]  {{\em #1\/} ({#2}, {#3} {#4})}
 \newcommand\inBO[8]{{\em #8}, in:\ {\em #1}, {#2}\ ({#3}, {#4} {#5}), p.\ {#6--#7}}
 \def\dim   {dimension}
 \def\jf    {J.\ Fuchs}
 \def\adma  {Adv.\wb Math.}
 \def\aspm  {Adv.\wb Stu\-dies\wB in\wB Pure\wB Math.}
 \def\aste  {Ast\'e\-ris\-que}
 \def\coma  {Con\-temp.\wb Math.}
 \def\Coma  {Con\-temp. Math.}
 \def\comp  {Com\-mun.\wb Math.\wb Phys.}
 \def\Comp  {Com\-mun.\wb Math. Phys.}
 \def\cOmp  {Com\-mun. Math.\wb Phys.}
 \def\cpma  {Com\-pos.\wb Math.}
 \def\Fiic  {Fields Institute\wB Commun.}
 \def\ijmp  {Int.\wb J.\wb Mod.\wb Phys.\ A}
 \def\ijMp  {Int.\wb J.\linebreak[0]Mod.\wb Phys.\ A}
 \def\joal  {J.\wB Al\-ge\-bra}
 \def\jofa  {J.\wb Funct.\wb Anal.}
 \def\joms  {J.\wb Math.\wb Sci.\wB Tokyo}
 \def\jpaa  {J.\wB Pure\wB Appl.\wb Alg.}
 \def\nuci  {Nuovo\wB Cim.}
 \def\nupb  {Nucl.\wb Phys.\ B}
 \def\Nupb  {Nucl. Phys.\ B}
 \def\phlb  {Phys.\wb Lett.\ B}
 \def\phrl  {Phys.\wb Rev.\wb Lett.}
 \def\plms  {Proc.\wB Lon\-don\wB Math.\wb Soc.}
 \def\pnas  {Proc.\wb Natl.\wb Acad.\wb Sci.\wb USA}
 \def\reth  {Re\-present.\wB Theory}
 \def\rims  {Publ.\wB RIMS}
 \def\sema  {Selecta\wB Mathematica}  
 \def\slnm  {Sprin\-ger\wB Lecture\wB Notes\wB in\wB Mathematics}
   \def\AMS    {{Ame\-ri\-can Mathematical Society}}
   \def\BIR    {{Birk\-h\"au\-ser}}
   \def\Bo     {{Boston}}
   \def\PR     {{Providence}}

\def\cats {categories}
\def\voa  {vertex operator algebra} 
\bibliographystyle{amsalpha}  \medskip
\end{document}